\documentclass[smallextended,envcountsect]{svjour3} 
\smartqed 
\usepackage{amsmath,amssymb}
\usepackage{array}
\usepackage{blkarray}
\usepackage{cite}
\usepackage{float}
\usepackage{graphicx}
\usepackage{hyperref}
\usepackage{pbox}
\usepackage{pdflscape}
\graphicspath{{Figures/}}
\floatstyle{boxed}
\newfloat{procedure}{thp}{lop}

\newtheorem{myAlgorithm}{Procedure}

\usepackage{xspace}

\newcommand{\head}[1]{\ensuremath{\mathtt{head}(#1)}\xspace}
\newcommand{\init}[1]{\ensuremath{\mathtt{init}(#1)}\xspace}
\newcommand{\mdeg}[1]{\ensuremath{\mathtt{mdeg}(#1)}\xspace}
\newcommand{\mvar}[1]{\ensuremath{\mathtt{mvar}(#1)}\xspace}
\newcommand{\mydeg}[2]{\ensuremath{\mathtt{deg}(#1,#2)}\xspace}
\newcommand{\rank}[1]{\ensuremath{\mathtt{rank}(#1)}\xspace}
\newcommand{\sat}[1]{\ensuremath{\mathtt{sat}(#1)}\xspace}
\newcommand{\tail}[1]{\ensuremath{\mathtt{tail}(#1)}\xspace}

\newcommand{\CalC}[1]{\ensuremath{\mathcal{C}^{#1}}\xspace}
\newcommand{\Man}{\ensuremath{\mathcal{M}}\xspace}
\newcommand{\Nash}{\ensuremath{\mathcal{N}}\xspace}
\newcommand{\Qz}{\ensuremath{\mathbb{Q}[\mathbf{z}]}\xspace}

\newcommand{\Rn}[1]{\ensuremath{\mathbb{R}^{#1}}\xspace}
\newcommand{\TMan}[1]{\ensuremath{\mathcal{T}_{#1}\mathcal{M}}\xspace}
\journalname{JOTA}

\begin{document}

\title{Whitney's Theorem, Triangular Sets and Probabilistic Descent on Manifolds}

\author{David W. Dreisigmeyer
%\thanks{Communicated by Alexandru Krist\'{a}ly.}
}

\institute{
  David W. Dreisigmeyer,  Corresponding author \at
  United States Census Bureau \\
  Center for Economic Studies \\
  Suitland, MD \\
  Department of Electrical and Computer Engineering \\
  Colorado State University \\
  Fort Collins, CO \\
  david.wayne.dreisigmeyer@census.gov
}

%The correct dates will be entered by the editor.

\maketitle

\begin{abstract}
We examine doing probabilistic descent over manifolds implicitly defined by a set of polynomials with rational coefficients.
The system of polynomials is assumed to be triangularized.
An application of Whitney's embedding theorem allows us to work in a reduced dimensional embedding space.
A numerical continuation method applied to the reduced-dimensional embedded manifold is used to drive the procedure.
\end{abstract}
\subclass{65K10 \and 90C56}

%============================================================
\section{\label{sec:introduction}Introduction}
%============================================================
In an optimization problem with equality constraints the feasible set often has nice geometric properties.
If we `stand' at a typical point the feasible set will locally look to us like a Euclidean space.
But this locally Euclidean space will be contained in a larger space.
For example, the feasible set could be a one-dimensional string contained in a hundred-dimensional ambient space.
It seems like a waste to work with a hundred variables when the geometric structure we're actually interested in is only one-dimensional.
Whiney's theorem tells us that in fact we could work in a three-dimensional ambient space.
The string can be projected from the original hundred dimensions into a random three-dimensional subspace.
Provided we can find a map from the projected string back to the original string we lose nothing by working in the three-dimensional space.

Sard's theorem is well known in optimization since it tells us that equality constraints often give a reasonable feasible set over which we work.
Or at least it will alert us when this may not be the case.
Whitney's theorem is perhaps less well known.
The theorem offers a potentially large decrease in the effective dimension one needs to optimize in, which is certainly an attractive proposal.
However this is not cost free.
Firstly, the decrease in dimensionality is achieved by a linear projection from the initial ambient space into a lower-dimensional working space.
This requires some sort of variable elimination from the set of equality constraints in order to do the projection.
The variable elimination is the nonlinear analogue of triangularizing a matrix.
Secondly, the reconstruction from the reduced dimension working space back into the original space is generally nonlinear.
A numerical implementation of the inverse (implicit) function theorem over a manifold is required.
This is a nonlinear version of basic and nonbasic variables.
The manifold in the reduced dimensional space serves as the basic variables and the embedding in the higher dimensional space represents the nonbasic variables.

A decade ago the author was interested in applying Whitney's theorem to optimization problems.
Unfortunately the required machinery was not quite sufficient enough to give any satisfactory results.
In the intervening years developments in polynomial decompositions over real numbers has made it possible to provide a sketch of the main ideas.
While the method is not quite complete it is complete enough to offer some potential applications and, more importantly, a clear map of what work still needs to be done.
In other words: \textit{it works but it may not yet work well enough}.

The paper is organized as follows.
We give a simple example in Section ~\ref{sec:motivate} that demonstrates the basic problem and principles we will be dealing with.
Section~\ref{sec:whitney} covers the necessary results from differential geometry.
The main result is Whitney's weak embedding theorem which gives an upper bound on the minimal embedding dimension of a manifold.
In Section~\ref{sec:triangular} we look at some needed results from real algebraic geometry.
In particular some real polynomials decompositions are introduced.
These put the constraint equations into a form convenient for applying Whitney's theorem, which we examine in Section~\ref{sec:tri-whit}.
Moving over manifolds in a controllable and intelligent way is necessary.
Section~\ref{sec:moving} looks at some methods to accomplish this.
The full optimization procedure is laid out in Section~\ref{sec:opt-method} with Section~\ref{sec:example} providing an example.
A discussion follows in Section~\ref{sec:discussion}.

%============================================================
\section{\label{sec:motivate}Motivating Example}
%============================================================
Here we give a simple example in order to motivate and demonstrate the methodology that will be developed.
Consider the optimization problem
\begin{align*}
\min_{\mathbf{z} \in \Rn{4}} & f(\mathbf{z})
\end{align*}
subject to the constraints
\begin{align*}
0 &= G_1(\mathbf{z}) = z_4 + z_3 + z_2 + z_1 \mbox,	\\
0 &= G_2(\mathbf{z}) = z_4z_1 + z_4z_3 + z_3z_2 + z_2z_1 \mbox,	\\
0 &= G_3(\mathbf{z}) = z_4z_2z_1 + z_4z_3z_1 + z_4z_3z_2 + z_3z_2z_1 \mbox{ and }	\\
0 &= G_4(\mathbf{z}) = z_4z_3z_2z_1 - 1 \mbox.
\end{align*}
The constraint set can be rewritten in the form \cite{art:kalkbrener-1993}
\begin{align}
\label{eq:ex_triangular}
\begin{split}
	0 &= g_1(x,u) = u^2 x^2 - 1	\\
	0 &= g_2(y_1,u) = y_1 + u	\\
	0 &= g_3(y_2,x) = y_2 + x
\end{split}
\end{align}
where we renamed the variables $u = z_1$, $x = z_2$, $y_1 = z_3$ and $y_2 = z_4$.
Looking at the polynomials in (\ref{eq:ex_triangular}) we see that they are in a `triangular form': $g_i$ only depends on $u$, $x$ or $y_j$ where $j < i$.

Looking at the new polynomials we see that it is not necessary to work in \Rn{4}.
Only the $(x,u) \in \Rn{2}$ need to be explicitly optimized over, treating $y_1(x,u) = - u$ and $y_2(x,u) = -x$ as functions of $u$ and $x$.
The original optimization problem now becomes
\begin{align*}
\min_{(x,u) \in \Rn{2}} \hat{f}(y_1(x,u),y_2(x,u),x,u)
\mbox{ subject to }
g_1(x,u) = 0
\end{align*}
where $f(z_1,z_2,z_3,z_4) = \hat{f}(z_3,z_4,z_2,z_1)$.
[We will typically drop the `hat' notation on the objective function when the variables are simply permuted.
The context should make it clear when this occurs.]

This example highlights the key techniques we will employ.
First, we have an optimization problem initially formulated in a `large' Euclidean space.
The constraints are given by polynomials with rational coefficients.
A triangular decomposition can be performed on the original constraint statements to derive a new set of constraints.
This new constraint set is then used to reduce the effective dimension we need to optimize in.
The optimization can now proceed in a lower-dimensional Euclidean space versus the original one.
As part of the projection into the lower-dimensional space an inverse function is used to map from the feasible set in the lower dimension back to the original feasible set.

%============================================================
\section{\label{sec:whitney}Whitney's Theorem}
%============================================================
A manifold \Man is a subset of \Rn{n} that looks locally like \Rn{m}, $n > m$.  
Let $\mathcal{I}$ be an index set and for each $i \in \mathcal{I}$ let $\mathcal{V}_{i} \subset \Man$ be an open subset where $\bigcup_{i \in \mathcal{I}} \mathcal{V}_{i} = \Man$.  
Additionally, let there be homeomorphisms $\phi_{i} : \mathcal{V}_{i} \rightarrow \Rn{m}$.  
Each pair $(\mathcal{V}_{i}, \phi_{i})$ is called a \emph{chart} and the set $\{ (\mathcal{V}_{i}, \phi_{i}) \}_{i \in \mathcal{I}}$ is called an \emph{atlas}.

When two open subsets $\mathcal{V}_{i}$ and $\mathcal{V}_{j}$ of \Man overlap there will be a homeomorphic transition function $\tau_{ij} : \Rn{m} \rightarrow \Rn{m}$ on the overlap $\mathcal{V}_{i} \cap \mathcal{V}_{j}$ defined by $\tau_{ij} = \phi_{j} \circ \phi_{i}^{-1}$.  
The transition functions allow us to stitch together the charts into a coherent whole.  
This then defines a \emph{topological manifold}.

The transition functions are often taken to be differentiable.  
In this case we have a \emph{differentiable manifold}.  
If the $\tau_{ij}$ are all \CalC{k}, $k \geq 1$, then \Man is said to be a \CalC{k}-manifold.  
If \Man is a \CalC{\infty}-manifold then it is called \emph{smooth}.  
Every \CalC{k}-manifold can be made smooth \cite{book:hirsch-1997}.

At every point $p \in \Man \subset \Rn{n}$ a tangent space \TMan{p} can be attached.  
If an inner product $h: \TMan{p} \times \TMan{p} \rightarrow \mathbb{R}$ is defined for every $p \in \Man$ and varies smoothly over \Man then we have a \emph{Riemannian manifold}.  
The function $h(\cdot,\cdot)$ is called a \emph{metric} and with it notions of lengths and angles in \TMan{p} are available over \Man.

The manifolds we consider arise from equality constraints.
\begin{definition}[Regular level sets]
\label{def:level-sets}
Let $\mathbf{g}: \Rn{m+k} \rightarrow \Rn{k}$.  
For the set $\Man = \{ \mathbf{z}:\mathbf{g}(\mathbf{z}) = \mathbf{c} \}$ assume $\nabla \mathbf{g}(\mathbf{z})$ is full rank.  
Then \Man is a \emph{regular level set} of $\mathbf{g}(\mathbf{z})$.  
Additionally, the null space of $\nabla \mathbf{g}(\mathbf{z})$ coincides with the tangent space to \Man.
$\square$
\end{definition}
This paper will consider equality constraints given by members of \Qz, the polynomials in $\mathbf{z} = [z_1,\ldots,z_{m+k}]$ with rational coefficients.
Then by Sard's theorem almost every level set is a Riemannian manifold.
\begin{theorem}[Regular level set manifolds and Sard's theorem \cite{book:hirsch-1997}]
\label{theo:sard}
Let $\mathbf{g}: \Rn{m+k} \rightarrow \Rn{k}$ be a \CalC{\infty} function.  
Then almost every level set of $\mathbf{g}(\mathbf{z})$ is regular and a $m$-dimensional manifold.
A regular level set acquires a metric from its embedding in \Rn{m+k} by restricting the standard Euclidean metric to the tangent space $\TMan{\mathbf{p}}$ given by $\mathrm{null}(\nabla \mathbf{g}(\mathbf{p}))$ for all $\mathbf{p} \in \Man$.
\end{theorem}

One issue is that the $k$ in Theorem~\ref{theo:sard} may be much larger than $m$.
So even though the manifold may be low-dimensional it could be embedded in a very high-dimensional space.
Whitney's theorem gives an upper bound on the minimal embedding dimensional for any manifold.
\begin{theorem}[Whitney's weak embedding theorem \cite{book:hirsch-1997}]
\label{theo:whitney}
Let \Man be an $m$-dimensional manifold.
Then \Man can be embedded in \Rn{2m+1}.
Further, if \Man is embedded in \Rn{m+k}, $k > m + 1$, a random projection into a $2m+1$ subspace gives an embedding of \Man in \Rn{2m+1}.
\end{theorem}

%============================================================
\section{\label{sec:triangular}Triangular Sets}
%============================================================
Here we cover the necessary background from real algebraic geometry.
Our main focus is putting a set of polynomials $P \subset \Qz$ into a convenient form for utilizing Whitney's theorem.
\begin{definition}[Semi-algebraic sets, real algebraic varieties and Nash manifolds]
\label{def:semialgebraic_sets}
A set $S \subset \Rn{n}$ is a \textit{semi-algebraic set} if it satisfies a set of polynomial equations $g_i(\mathbf{z}) = 0$, $i=1,\ldots,j$, and polynomial strict inequalities $h_i(\mathbf{z}) > 0$, $i=1,\ldots,l$, or a finite union of such sets.
% This exhausts the relationships since we can replace inequalities and inequations by 
% \begin{align*}
% g(\mathbf{z}) \geq 0 &\rightarrow g(\mathbf{z}) - y^2 = 0 \mbox{, $y \in \Rn{}$, and }	\\
% h(\mathbf{z}) \neq 0 &\rightarrow y h(\mathbf{z}) > 0 \mbox{, $y \in \Rn{}$.}
% \end{align*}
\noindent A \textit{real algebraic variety} is defined using only a set of polynomial equations $g_i(\mathbf{z}) = 0$, $i=1,\ldots,j$, without finite unions.
\textit{Nash manifolds} are semi-algebraic sets that are also manifolds.
$\square$
\end{definition}
This paper focuses on Nash manifolds that are also varieties.

In order to utilize Whitney's theorem the polynomials defining a Nash manifold will need to be decomposed into a friendlier form.
Let the Nash manifold $\Nash \subset \Rn{m+k}$ be defined by the set of polynomial equations $g_i(\mathbf{z}) = 0$, $i=1,\ldots,k$.
Assign to the entries of $\mathbf{z}$ the ordering 
\begin{align*}
z_1 < z_2 < \ldots < z_{m+k} \mbox.
\end{align*}
Define $Q_i = \mathbb{Q}[z_1,\ldots,z_i]$ as the ring of polynomials in the variables $z_1,\ldots,z_i$ with rational coefficients.
\mydeg{g}{z_i} is the \textit{degree} of polynomial $g$ in the variable $z_i$.
\mvar{g} is the \textit{greatest variable} $z_i$ such that $\mydeg{g}{z_i} \neq 0$.
We can decompose a $g \in Q_i$ into the form $g = \iota z_i^d + \tau$ where $d = \mydeg{g}{z_i} > 0$, $\iota \in Q_{i-1}$ and $\tau \in Q_i$ with $\mydeg{\tau}{z_i} < d$.
Then $\init{g} = \iota$ is the \textit{initial} of $g$, $\mdeg{g} = d$ is the \textit{main degree} of $g$, $\rank{g} = z_i^d$ is the \textit{rank} of $g$, $\tail{g} = \tau$ is the \textit{tail} of $g$ and $\head{g} = g - \tail{g} = \iota z_i^d$ is the \textit{head} of $g$.
Every polynomial is then decomposed as $p = \init{p} \rank{p} + \tail{p}$.
\begin{definition}[Triangular sets]
\label{def:triangular_sets}
Let $T \subset Q_n$.
We call $T$ a \textit{triangular set} if no member of $T$ is constant and for any pair $g,h \in T$, $g \neq h$, we have $\mvar{g} \neq \mvar{h}$.
If $z_i = \mvar{g}$ for some $g \in T$ then $z_i$ is called \textit{algebraic with respect to} $T$.
Otherwise $z_i$ is called \textit{free with respect to} $T$.
$\square$
\end{definition}
Triangular sets are a nonlinear generalization of triangulating a matrix.
They are very convenient for using Whitney's theorem.
\begin{example}
Let $T = \{ z_1^2 z_2^2 -1, z_3 + z_1, z_4 + z_2 \}$.
This is a triangular set where 
\begin{align*}
\mvar{z_1^2 z_2^2 -1} &= z_2 \mbox, \\
\mvar{z_3 + z_1} &= z_3 \mbox{ and }\\
\mvar{z_4 + z_2} &= z_4 \mbox.
\end{align*}
Then $z_1$ is free and the other variables are algebraic.
$\square$
\end{example}

One of the difficulties we face with using triangular sets is that most of the prior work on them has been with $\mathbb{C}[\mathbf{z}]$, the polynomials with complex coefficients.
Our concern is only with the real solutions of a system of polynomials in \Qz which has seen recent development \cite{art:chen-2011, art:chen-2013}.
We'll develop the necessary machinery for working with Nash manifolds but first we'll look at a standard solution method over $\mathbb{C}$.

The algorithm developed by Kalkbrener in \cite{art:kalkbrener-1993} seems like a reasonable choice for calculating triangular sets \cite{art:aubry-1999b,art:chen-2012}.
This method builds on a particular type of triangular set which requires some additional machinery.
\begin{definition}[Saturated ideals]
\label{def:sat_ideals}
Given a triangular set $T = \{t_1,\ldots,t_k\}$ the \textit{ideal} $\langle T \rangle$ generated by $T$ is given by
\begin{align*}
\langle T \rangle &=
	\left\lbrace
    	q \in \mathbb{Q}[\mathbf{z}] :
        	q = \textstyle\sum_i p_i t_i \mbox{, } p_i \in \mathbb{Q}[\mathbf{z}]
    \right\rbrace \mbox.
\end{align*}
Let $h_T = \textstyle\prod_i \init{t_i}$.
The \textit{saturated ideal} of $T$ is defined as
\begin{align*}
\sat{T} &=
	\left\lbrace
    	q \in \mathbb{Q}[\mathbf{z}] :
        	\exists n \in \mathbb{N}_0 \mbox{ such that } h_T^n q \in \langle T \rangle
    \right\rbrace
\end{align*}
with $\sat{\varnothing} = \langle 0 \rangle$.
$\square$
\end{definition}
The ideal $\langle T \rangle$ gives all valid equality constraints that follow from those in $T$.
We see that $\langle T \rangle \subset \sat{T}$ but \sat{T} can contain additional polynomials.
% The question of when $\langle T \rangle = \sat{T}$ has been examined in \cite{art:lemaire-2008}.
% \sat{T} is of greater relevance to Kalkbrener's method than $\langle T \rangle$.
% From this it follows that the variety defined by $T$ must contain the variety defined by \sat{T}.
% The differences between the two varieties is important for defining well-behaved triangular sets.

\begin{definition}[Regular chains]
\label{def:regular_chain}
Let $T = \{t_1,\ldots,t_k\}$ be a triangular set.
A $p \in \Qz$ is \textit{regular modulo} \sat{T} if $p \not\in \sat{T}$ and there does not exist a $q \in \Qz$ where $q \not\in \sat{T}$ but $pq \in \sat{T}$.
Define $T_{<j} = \{t_1,\ldots,t_{j-1}\}$, $j < k$.
$T$ is a \textit{regular chain} if $T = \varnothing$, or $T_{<k}$ is a regular chain and \init{t_k} is regular modulo \sat{T_{<k}}.
$\square$
\end{definition}
\begin{definition}[Kalkbrener triangular decomposition\cite{art:chen-2013}]
\label{def:kalkbrener}
Let $\mathcal{T} = \{T_1,\ldots,T_l\}$ be a set of regular chains and $V(F)$ be an algebraic variety.
$\mathcal{T}$ is a \textit{Kalkbrener triangular decomposition} of $V(F)$ if $V(F) = \bigcup_i V(\sat{T_i})$.
$\square$
\end{definition}
So $V(\sat{T_i})$ is geometrically more relevant to Kalkbrener's method than $V(T_i)$: each $V(\sat{T_i})$ gives some piece of $V(F)$.
We can relate $V(T)$ and $V(\sat{T})$.
\begin{definition}[Regular zeros]
\label{def:regular_zeros}
Given a triangular set $T = \{t_1,\ldots,t_k\}$ let $h_T = \textstyle\prod_i \init{t_i}$.
The \textit{regular zeros} $W(T)$ of the triangular set $T$ is defined as $W(T) = V(T) - V(h_T)$.
$\square$
\end{definition}
For a triangular set $T$ it is the case that \cite{art:aubry-1999a}
\begin{align*}
V(\sat{T}) &= \overline{W(T)}
\end{align*}
where closure is with respect to the Zariski topology.
If $V(T)$ is also a manifold then $\overline{W(T)} = V(T) = V(\sat{T})$.
\begin{definition}[Radical of $\langle T \rangle$]
\label{def:radical}
The \textit{radical} of the ideal $\langle T \rangle$ is
\begin{align*}
\sqrt{\langle T \rangle} &= 
	\left\lbrace q \in \Qz :
		\exists n \in \mathbb{N}_0 \mbox{ such that } q^n \in \langle T \rangle
    \right\rbrace \mbox.
\end{align*}
That is, $\sqrt{\langle T \rangle}$ consists of all those polynomials that vanish on $V(T)$.
$\square$
\end{definition}
Since $V(T) = V(\sat{T})$ it follows from Hilbert's Nullstellensatz \cite{book:cox-2015} that $\sqrt{\langle T \rangle} = \sqrt{\sat{T}}$.
From this we see that a Kalkbrener triangular decomposition $\mathcal{T} = \{T_1,\ldots,T_l\}$ of $F \subset \Qz$ gives $\sqrt{\langle F \rangle} = \bigcap_i \sqrt{\sat{T_i}}$.

A Kalkbrener triangular decomposition allows for complex solutions but Nash manifolds only work with the real solutions.

%============================================================
\section{\label{sec:tri-whit}Whitney's Theorem Applied to Triangular Sets}
%============================================================
We first look at an easy to understand case.
Consider the linear system $A\mathbf{z} = \mathbf{b}$.
This set of polynomials defines a special type of manifold that is globally like \Rn{m} when $A \in \Rn{k \times (m+k)}$.
Let $k > m+1$ for what follows.
It is convenient to have $A$ in a triangular form
\begin{align*}
A &= 
    \begin{blockarray}{cccc}
    	k-m-1 & m+1 & m & \\
    	\begin{block}{[ccc]c}
          A_{11} & A_{12} & A_{13} & k-m-1 \\
          0 & A_{22} & A_{23} & m+1 \\
    	\end{block}
    \end{blockarray}
\end{align*}
with $A_{11}$ and $A_{22}$ upper triangular.
Partition the $\mathbf{z}$ as
\begin{align*}
\mathbf{z} &=[\mathbf{y}\ \mathbf{x}\ \mathbf{u}]
\end{align*}
with $\mathbf{y} \in \Rn{k-m-1}$, $\mathbf{x} \in \Rn{m+1}$ and $\mathbf{u} \in \Rn{m}$.
Whitney's theorem tells us that the space given by $\mathbf{x}$ and $\mathbf{u}$ is sufficient to embed the manifold implicitly defined by the set of polynomials $A\mathbf{z} - \mathbf{b}$.
Let the original manifold defined in \Rn{m+k} be denoted by \Man and the projected image in \Rn{2m+1} be denoted by $\Man^{\prime}$.
Partition $\mathbf{b}$ as $\mathbf{b} = [\mathbf{b}_1\ \mathbf{b}_2]$.
We can solve the equation
\begin{align*}
	\left[A_{22}\ A_{23}\right]
      \left[\begin{array}{c}
       \mathbf{x} \\
       \mathbf{u}
      \end{array}\right]
	&=
    \mathbf{b}_2
\end{align*}
to find the points on $\Man^{\prime}$.
Call one of these $[\mathbf{x}^*\ \mathbf{u}^*]^T \in \Man^{\prime}$. 
This can then be mapped onto the corresponding point on \Man by solving
\begin{align*}
	\left[A_{11}\ A_{12}\ A_{13}\right]
    \left[\begin{array}{c}
      \mathbf{y} \\
      \mathbf{x}^* \\
      \mathbf{u}^*
    \end{array}\right]
	&=
    \mathbf{b}_1
\end{align*}
for $\mathbf{y}$.
While we may not want to do it in this case, the optimization can be carried out in \Rn{2m+1} and we can then map from $\Man^{\prime}$ to \Man by finding $\mathbf{y} = \mathbf{y}(\mathbf{x},\mathbf{u})$.
Using the full vector $\mathbf{z} = [\mathbf{y}\ \mathbf{x}\ \mathbf{u}]$ we can pull back any function defined on \Man to $\Man^{\prime}$.
A function that is pulled back could be, e.g., the objective function or an inequality constraint.

Now we can repeat the above example for the general nonlinear case.
We will make the simplifying assumption that \Man is path connected.
If there are multiple pieces of \Man then the following argument still carries through for each individual piece.
After a possible reordering of the coordinates, partition the $\mathbf{z} \in \Rn{m+k}$ from Definition~\ref{def:level-sets} as 
\begin{align*}
\mathbf{z} &=[\mathbf{y}\ \mathbf{x}\ \mathbf{u}]
\end{align*}
where $\mathbf{y} \in \Rn{k-m-1}$, $\mathbf{x} \in \Rn{m+1}$ and $\mathbf{u} \in \Rn{m}$.
We take the members of $\mathbf{g}(\mathbf{z})$ to be triangularized \cite{art:chen-2013} in the $\mathbf{x}$ and $\mathbf{y}$ as
\begin{subequations}
\label{eq:g_star_and_ast}
\begin{align}
\mathbf{g}^{\star}(\mathbf{x}, \mathbf{u}) : \Rn{2m+1} &\rightarrow \Rn{m+1}
	\mbox { where }			
    \label{eq:g_star-def} \\
g^{\star}_i(\mathbf{x}, \mathbf{u}) &= g^{\star}_i(x_1,\ldots,x_i,\mathbf{u})
	\label{eq:g_star_i-def}
\end{align}
for $i=1,\ldots,m+1$ and
\begin{align}
\mathbf{g}^{\circ}(\mathbf{y};\mathbf{x}, \mathbf{u}) : \Rn{k-m-1} &\rightarrow \Rn{k-m-1} 
	\mbox { where }			
    \label{eq:g_ast-def} \\
g^{\circ}_j(\mathbf{y};\mathbf{x}, \mathbf{u}) &= g^{\circ}_j(y_1,\ldots,y_j;\mathbf{x}, \mathbf{u})
	\label{eq:g_ast_j-def}
\end{align}
\end{subequations}
for $j=1,\ldots,k-m-1$, with $\mathbf{g} = [\mathbf{g}^{\star}\ \mathbf{g}^{\circ}]^T$.
In (\ref{eq:g_ast-def}) the $\mathbf{x}$ and $\mathbf{u}$ are determined from (\ref{eq:g_star-def}) and are treated as parameters.
If there is at most one solution to $\mathbf{g}^{\circ}(\mathbf{y}; \mathbf{x}, \mathbf{u}) = \mathbf{0}$ for every $\mathbf{x}$ and $\mathbf{u}$ such that $\mathbf{g}^{\star}(\mathbf{x}, \mathbf{u}) = \mathbf{0}$, then $\mathbf{g}^{\star}(\mathbf{x}, \mathbf{u})$ gives a low-dimensional embedding of our manifold.
Theorem~\ref{theo:whitney} tells us that typically any choice for the $\mathbf{x}$ and $\mathbf{u}$ should accomplish this.

We see that the optimization can progress in \Rn{2m+1} by using the $\mathbf{x}$ and $\mathbf{u}$ variables.
The manifold defined implicitly by (\ref{eq:g_star-def}) can then be mapped onto the original manifold by solving (\ref{eq:g_ast-def}) for the unique $\mathbf{y}(\mathbf{x},\mathbf{u})$.
We then use
\begin{align*}
\mathbf{z}(\mathbf{u},\mathbf{x}) &= 
	[\mathbf{y}(\mathbf{x},\mathbf{u})\ \mathbf{x}\ \mathbf{u}]
\end{align*}
to pull back the objective function and any inequality constraints from the manifold embedding in \Rn{m+k} to the embedding in \Rn{2m+1}.

The ability to triangularize the polynomials in $\mathbf{g}(\mathbf{z})$ is what allows Whitney's theorem to be easily implemented.
There are two things to notice about this triangularization.
First, all of this is restricted to working only with the real versus complex roots.
Second, as noted above, we only work on the highest dimensional subcomponent of the variety $V(\mathbf{g})$.
% In \cite{art:chen2013} this is called a lazy triangular decomposition.

Finding the function $\mathbf{y}(\mathbf{x},\mathbf{u})$ is the crucial step.
The most obvious way of doing this is the one we presented above.
However it's unimportant how this is done in practice.
Any numerical implementation of the implicit function theorem applied to $\mathbf{g}^{\circ}(\mathbf{y}, \mathbf{x}, \mathbf{u})$ will suffice, where now $\mathbf{x}$ and $\mathbf{u}$ are treated as variables.

%============================================================
\section{\label{sec:moving}Moving Over the Manifold}
%============================================================
We will need a procedure to intelligently move over the manifold $\Man^{\prime} \subset \Rn{2m+1}$.
One advantage of working with a lower-dimensional embedding space is that methods that may be expensive when working in high-dimensions can become feasible.
This is counter-balanced by the need to map back into the original space, say with the $\mathbf{y}(\mathbf{x},\mathbf{u})$ from Section~\ref{sec:triangular}.
But we're free to use different methods to accomplish these different goals.
We move along a path on the manifold $\Man^{\prime}$ with one technique and drag along a point on $\Man \subset \Rn{m+k}$ with a different one.
Numerical continuation methods will likely play a key roll in any implementation employing Whitney's theorem.
We mention the references \cite{art:allgower-1985,art:allgower-2000,art:rheinboldt-2000} in addition to the ones below.

On a Riemannian manifold \Man a \emph{geodesic} is the shortest path between two (nearby) points on \Man and generalizes the notion of a straight line.
The coupled nonlinear equations for a geodesic on \Man are \cite{misc:dreisigmeyer-2006}
\begin{align*}
% \label{eq:geodesics}
\ddot{z}^i + \Gamma^{i}_{jk} \dot{z}^j \dot{z}^k &= 0 \\
\Gamma^{i}_{jk} &=
	\left[ \nabla \mathbf{g}^+ \right]^i 
    \frac{\partial^2 \mathbf{g}}{\partial z^j \partial z^k}
\end{align*}
where $\left[ \nabla \mathbf{g}^+ \right]^i$ is the $i$-th row of the pseudo-inverse of $\nabla \mathbf{g}$.
The Einstein summation convention is used above where a repeated index in a term is summed over.
The $\Gamma^{i}_{jk}$ are the \textit{Christoffel symbols} associated with $\mathbf{g}$.
The distance moved over the manifold in a unit time period is controlled by the length of the initial tangent vector $\dot{\mathbf{z}}_0$.
For $\Man^{\prime} \subset \Rn{2m+1}$ $\mathbf{g}^\star$ is substituted for $\mathbf{g}$ and we only work with the $\mathbf{x}$ and $\mathbf{u}$ variables.
If $2m+1$ is `small' it's possible that finding geodesics on $\Man^{\prime}$ would be feasible.
One could then find $\mathbf{y}(\mathbf{x},\mathbf{u})$ for the end point of a geodesic and use this to pull back any necessary functions from \Man to $\Man^{\prime}$.

In \cite{art:brodzik-1998} a method of approximating a manifold by locally projecting down tangent spaces was developed.
% (See \cite{art:rheinboldt-2000} also.)
Very little of the machinery is needed.
As stated in Theorem~\ref{theo:sard} the tangent space $\TMan{\mathbf{p}}^\prime$ to $\Man^\prime$ at $\mathbf{p} = [\mathbf{x}\ \mathbf{u}] \in \Man^\prime$ is given by $\mathrm{null}(\nabla \mathbf{g}^\star(\mathbf{p}))$.
Let $U^\prime$ be an orthonormal basis for $\TMan{\mathbf{p}}^\prime$.
Then the vector $\mathbf{q}_0 = \mathbf{p} + U^\prime \mathbf{w}$ lies in $\TMan{\mathbf{p}}^\prime$ when the latter is viewed as an affine subspace of \Rn{2m+1}.
Set $N^\prime = [\nabla \mathbf{g}^\star(\mathbf{p})]^+$ where $A^+$ is the Moore-Penrose pseudoinverse of the matrix $A$.
The tangent vector $\mathbf{q}_0$ is projected down onto $\Man^\prime$ by
\begin{subequations}
\label{eq:project-tan}
\begin{align}
\mathbf{q}_{n+1} &= \mathbf{q}_n - N^\prime \mathbf{g}^\star(\mathbf{q}_n) 
	\mbox{, $n=0,1,\ldots$, with}
	\label{eq:xn-def} \\
\mathbf{q}_0 &= \mathbf{p} + U^\prime \mathbf{w}
	\label{eq:x0-def}
\end{align}
\end{subequations}
for some $\mathbf{w} \in \Rn{m}$.
Computationally this is the most attractive of the procedures presented and is the one we will use for the remainder.
Since this is a numerical implementation of the implicit function theorem there will be a radius around the origin in $\TMan{\mathbf{p}}^\prime$ where it will be guaranteed to converge to a point on $\Man^\prime$.
This is restated as a restriction on the length of the $\mathbf{w}$ in (\ref{eq:x0-def}).

%============================================================
\section{\label{sec:opt-method}The Optimization Procedure}
%============================================================
All of the pieces are in place to state the optimization procedure.
Let the optimization problem be stated as
\begin{equation}
\label{prob:gen-opt}
    \min_{\mathbf{z} \in \Rn{m+k}} \ f(\mathbf{z})
    \mbox{ subject to } \tilde{\mathbf{g}}(\mathbf{z}) \leq \mathbf{0}
\end{equation}
with $\tilde{\mathbf{g}}:\Rn{m+k} \rightarrow \Rn{k}$ and $\tilde{g}_i \in \mathbb{Q}[\mathbf{z}]$.
Triangularize $\tilde{\mathbf{g}}(\mathbf{z})$ into $\mathbf{g}(\mathbf{z})$ which has the structure presented in (\ref{eq:g_star_and_ast}).
Take $\mathbf{g}^{\star}: \Rn{2m+1} \rightarrow \Rn{m+1}$ as defining the manifold $\Man^\prime$.

We will assume we have a black-box implementation $\mathcal{B}(U^\prime,\mathbf{w})$ of the mapping defined by (\ref{eq:project-tan}).
$\mathcal{B}(U^\prime,\mathbf{w})$ will return \texttt{FALSE} if the algorithm doesn't converge to a point on $\Man^\prime$ close to the base point $\mathbf{p} \in \Man^\prime$.
That is, $\mathcal{B}(U^\prime,\mathbf{w})$ is also an oracle for when the implicit function theorem no longer holds around $\mathbf{p}$.
We will additionally assume a black-box implementation $\Phi:\Man^\prime \rightarrow \Man$ of the implicit function mapping.
For $\mathbf{p} \in \Man^\prime$ the point $\Phi(\mathbf{p}) = [\mathbf{y}(\mathbf{p})\ \mathbf{p}]$ is the unique point that lies on \Man.

Now we will use the probabilistic descent method developed in \cite{art:gratton-2015}.
A particular implementation is given by Procedure~\ref{alg:final_alg}.
The \textbf{Move} step in Procedure~\ref{alg:final_alg} picks a new base point on $\Man^\prime$ when the oracle $\mathcal{B}(U^\prime,\mathbf{w})$ tells us that the implicit function implementation is failing.
The base point defines the tangent space we are working in.
If the base point eventually becomes fixed then the convergence results are the same as those presented in \cite{art:gratton-2015}.
If the base point never becomes fixed we say the algorithm didn't converge and no solution was found.

\begin{procedure}[hbt!]
    \begin{myAlgorithm}[Probabilistic descent on manifolds]
    \label{alg:final_alg}
Specify the maximum step size increase $\alpha_{max} \in (0,\infty]$ and an $\alpha_0 \in (0,\alpha_{max}]$.
Set $\theta = 1/2$ and $\gamma = 2$.
Let the forcing function be $\rho(\alpha) = C \alpha^2$ for some constant $C$.
Start with an initial point $\mathbf{p}_0 \in \Man^\prime$ and set the current base point $\mathbf{b}_0 = \mathbf{p}_0$ and the current tangent vector $\mathbf{w}_0 = \mathbf{0} \in \Rn{m}$.
Let $\tilde{f} = f \circ \Phi (\mathbf{p}) : \Man^\prime \rightarrow \Man$.
Pick a maximum number of iterations $j_{max} \in \mathbb{N}$.
Finally let $j = 0$.
      \begin{description}[font=\bfseries]
      	\item [Tangents]
        	Let $U_j^\prime$ be an orthonormal basis for $\TMan{\mathbf{b}_j}^\prime$.
        	Pick a random $\mathbf{u} \in \Rn{m}$ such that $\|\mathbf{u}\|=1$.
        	Construct the two vectors $\mathbf{w}_{\pm} = \mathbf{w}_j \pm \alpha_j \mathbf{u}$.
        \item [Move]
        	If $\mathcal{B}(U_j^\prime,\mathbf{w}_+)$ or $\mathcal{B}(U_j^\prime,\mathbf{w}_-)$ returns \texttt{FALSE} set $\mathbf{p}_{j+1} = \mathbf{p}_j$, $\mathbf{b}_{j+1} = \mathbf{p}_j$, $\mathbf{w}_{j+1} = \mathbf{0}$ and $\alpha_{j+1} = \theta \alpha_j$.
            If $j > j_{max}$ return $\mathbf{p}_j$ otherwise increment $j$ and Goto the \textbf{Tangents} step.
		\item [Polling]
        	Let $\mathbf{p}_+ = \mathcal{B}(U_j,\mathbf{w}_+)$ and $\mathbf{p}_- = \mathcal{B}(U_j,\mathbf{w}_-)$.  
            If $\tilde{f}(\mathbf{p}_+) < \tilde{f}(\mathbf{p}_j) - \rho(\alpha_j)$ then set $\mathbf{p}_{j+1} = \mathbf{p}_+$, $\mathbf{w}_j = \mathbf{w}_+$ and declare the poll successful.  
            Else if $\tilde{f}(\mathbf{p}_-) < \tilde{f}(\mathbf{p}_j) - \rho(\alpha_j)$ then set $\mathbf{p}_{j+1} = \mathbf{p}_-$, , $\mathbf{w}_j = \mathbf{w}_-$ and declare the poll successful.  
            Otherwise set $\mathbf{p}_{j+1} = \mathbf{p}_j$ and declare the poll unsuccessful.
        \item [Step Size] 
        	If the \textbf{Polling} step was successful set $\alpha_{j+1} = \min(\alpha_{max},\gamma \alpha_j)$.
            Otherwise set $\alpha_{j+1} = \theta \alpha_j$.
            Set $\mathbf{b}_j = \mathbf{b}_{j-1}$.
            If $j > j_{max}$ return $\mathbf{p}_j$ otherwise increment $j$ and Goto the \textbf{Tangents} step.
      \end{description} 
    \end{myAlgorithm}
\end{procedure}

Because we are working with the vector space $\TMan{\mathbf{p}}^\prime$ we can also lay out a mesh.
This means that the MADS algorithm first developed in \cite{art:audet-2006} is also viable.
The advantage here is the ability to handle inequality constraints as well.

%============================================================
\section{\label{sec:example}An Example}
\begin{figure}[t]
	\centering
	\includegraphics[width=0.6\textwidth]{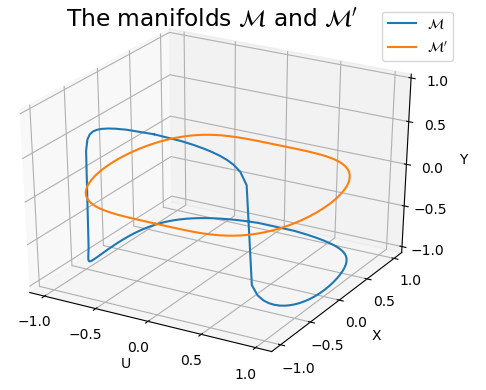}
  \caption{
  	The manifolds \Man and $\Man^\prime$ for the example in Section~\ref{sec:example}.
    $\Man^\prime$ is only embedded in \Rn{2} given by the coordinates $[x\ u]$.}
  \label{fig:example-M}
\end{figure}

%============================================================
Here we look at a simple example that demonstrates how to find the pieces for Procedure~\ref{alg:final_alg}.
Let the equality constraints be given by
\begin{align*}
\mathbf{g}(y,x,u) &= 
	\left[
    	\begin{array}{c}
	    	u^2 + x^3 + y^5		\\
    		u^4 + x^2 - 1
    	\end{array}
    \right] \\
    &= \mathbf{0}.
\end{align*}
See Figure~\ref{fig:example-M} for the solution set of $\mathbf{g}(y,x,u)$.
This triangularized system of polynomials gives the pair
\begin{align*}
g^\star(x,u) &= u^4 + x^2 - 1 \mbox{ and }\\
g^\circ(y;x,u) &= u^2 + x^3 + y^5 \mbox.
\end{align*}
Using $g^\circ(y;x,u)$ we can find $y$ as a function of $x$ and $u$:
\begin{align*}
y(x,u) &= -\sqrt[5]{u^2 + x^3} \mbox.
\end{align*}
This is the implicit function we need to map from $\Man^\prime$ back to \Man.
From this we have $\Phi = [y(x,u)\ x\ u]$.
We also need the gradient of $g^\star(x,u)$ to find the tangent space of $\Man^\prime$:
\begin{align*}
\nabla g^\star(x,u) &=
	\left[
    	\begin{array}{c}
            2 x		\\
    		4 u^3
    	\end{array}
    \right] \mbox.
\end{align*}
Looking at Figure~\ref{fig:example-M} we see that at least two different tangent spaces are required to implement the implicit function mapping $\TMan{\mathbf{p}}^\prime$ to $\Man^\prime$.
For $\mathcal{B}(U^\prime,\mathbf{w})$ we could have it return \texttt{FALSE} whenever, say, $\| \mathbf{q}_n - \mathbf{q}_0 \| > 1/2$, with the $\mathbf{q}_i$ as in (\ref{eq:project-tan}).
This would require more than two different tangent spaces to cover $\Man^\prime$ but it's likely an acceptable criterion to function as an oracle.

%============================================================
\section{\label{sec:discussion}Conclusions}
%============================================================
We have not implemented the procedure in Section~\ref{sec:opt-method} numerically.
There are multiple components to the method each of which is a significant undertaking in itself to be done correctly and efficiently.
However this was not our main purpose.
Rather we wished to draw upon classic and recent results that work together in ways that are perhaps novel for nonlinear optimization problems.

\bigskip
\footnotesize
\textbf{Disclaimer}
Any opinions and conclusions expressed herein are those of the author and do not necessarily represent the views of the U.S. Census Bureau.  
The research in this paper does not use any confidential Census Bureau information.  
This was authored by an employee of the US national government.  
As such, the Government retains a nonexclusive, royalty-free right to publish or reproduce this article, or to allow others to do so, for Government purposes only.

% \textbf{Acknowledgements} 
% The author would like to thank the editor and reviewers for their helpful comments and suggestions.
\normalsize

% \theendnotes
\bibliographystyle{spmpsci_unsrt}
\bibliography{jota_rhs}

\begin{thebibliography}{10}
\providecommand{\url}[1]{{#1}}
\providecommand{\urlprefix}{URL }
\expandafter\ifx\csname urlstyle\endcsname\relax
  \providecommand{\doi}[1]{DOI~\discretionary{}{}{}#1}\else
  \providecommand{\doi}{DOI~\discretionary{}{}{}\begingroup
  \urlstyle{rm}\Url}\fi

\bibitem{art:kalkbrener-1993}
Kalkbrener, M.: A generalized euclidean algorithm for computing triangular
  representations of algebraic varieties.
\newblock Journal of Symbolic Computation \textbf{15}(2), 143 -- 167 (1993)

\bibitem{book:hirsch-1997}
Hirsch, M.W.: Differential Topology.
\newblock Graduate Texts in Mathematics. Springer New York (1997)

\bibitem{art:chen-2011}
Chen, C., Davenport, J.H., Moreno~Maza, M., Xia, B., Xiao, R.: Computing with
  semi-algebraic sets represented by triangular decomposition.
\newblock In: Proceedings of the 36th International Symposium on Symbolic and
  Algebraic Computation, pp. 75--82 (2011)

\bibitem{art:chen-2013}
Chen, C., Davenport, J.H., May, J.P., Maza, M.M., Xia, B., Xiao, R.: Triangular
  decomposition of semi-algebraic systems.
\newblock Journal of Symbolic Computation \textbf{49}, 3 -- 26 (2013)

\bibitem{art:aubry-1999b}
Aubry, P., Maza, M.M.: Triangular sets for solving polynomial systems: a
  comparative implementation of four methods.
\newblock Journal of Symbolic Computation \textbf{28}(1), 125 -- 154 (1999)

\bibitem{art:chen-2012}
Chen, C., Maza, M.M.: Algorithms for computing triangular decomposition of
  polynomial systems.
\newblock Journal of Symbolic Computation \textbf{47}(6), 610 -- 642 (2012)

\bibitem{art:aubry-1999a}
Aubry, P., Lazard, D., Maza, M.M.: On the theories of triangular sets.
\newblock Journal of Symbolic Computation \textbf{28}(1), 105--124 (1999)

\bibitem{book:cox-2015}
Cox, D.A., Little, J., O'Shea, D.: Ideals, Varieties, and Algorithms, 4th edn.
\newblock Undergraduate Texts in Mathematics. Springer International Publishing
  (2015)

\bibitem{art:allgower-1985}
Allgower, E.L., Schmidt, P.H.: An algorithm for piecewise-linear approximation
  of an implicitly defined manifold.
\newblock SIAM Journal on Numerical Analysis \textbf{22}(2), 322--346 (1985)

\bibitem{art:allgower-2000}
Allgower, E.L., Georg, K.: Piecewise linear methods for nonlinear equations and
  optimization.
\newblock Journal of Computational and Applied Mathematics \textbf{124}(1), 245
  -- 261 (2000)

\bibitem{art:rheinboldt-2000}
Rheinboldt, W.C.: Numerical continuation methods: a perspective.
\newblock Journal of Computational and Applied Mathematics \textbf{124}(1), 229
  -- 244 (2000)

\bibitem{misc:dreisigmeyer-2006}
Dreisigmeyer, D.W.: Equality constraints, {R}iemannian manifolds and direct
  search methods.
\newblock
  \texttt{http://www.optimization-online.org/DB\_HTML/2007/08/1743.html} (2006)

\bibitem{art:brodzik-1998}
Brodzik, M.: The computation of simplicial approximations of implicitly defined
  p-dimensional manifolds.
\newblock Computers \& Mathematics with Applications \textbf{36}(6), 93 -- 113
  (1998)

\bibitem{art:gratton-2015}
Gratton, S., Royer, C.W., Vicente, L.N., Zhang, Z.: Direct search based on
  probabilistic descent.
\newblock SIAM Journal on Optimization \textbf{25}(3), 1515--1541 (2015)

\bibitem{art:audet-2006}
Audet, C., J.~E.~Dennis, J.: Mesh adaptive direct search algorithms for
  constrained optimization.
\newblock SIAM Journal on Optimization \textbf{17}(1), 188--217 (2006)

\end{thebibliography}

\end{document}